# Robust optimization Design of a New Combined Median Barrier Based on Taguchi method and Grey Relational Analysis

Yupeng Huang [a,b], Song Yao [a], Peng Chen [b], Zhengbao Lei [b,*], Xinzhong Gan [c], Youwei Gan [c]

[1] Key Laboratory of Traffic Safety on Track, Ministry of Education, School of Traffic and Transportation Engineering, Central South University, Hunan Changsha, 410075, China

[2] Large-scale Structural Impact Laboratory, Changsha University of Science and Technology, Hunan Changsha, 410114, China

[3] Yichun Highway Administration, Jiangxi Yichun, 336000, China

[4] Jiangxi Transportation Design and Research Institute Co., Ltd., Jiangxi Nanchang, 330052, China

**Abstract**：Accidents that vehicles cross median and enter opposite lane happen frequently, and the existing median barrier has weak anti-collision strength. A new combined median barrier (NCMB) consisted of W-beam guardrail and concrete structure was proposed to decrease deformation and enhance anti-collision strength in this paper. However, there were some uncertainties in the initial design of the NCMB. If the uncertainties were not considered in the design process, the optimization objectives were especially sensitive to the small fluctuation of the variables, and it might result in design failure. For this purpose, the acceleration and deflection were taken as objectives; post thickness, W-beam thickness and post spacing were chosen as design variables; the velocity, mass of vehicle and the yield stress of barrier components were taken as noise factors, a multi-objective robust optimization is carried out for the NCMB based on Taguchi and grey relational analysis (GRA). The results indicate that the acceleration and deflection after optimization are reduced by 47.3% and 76.7% respectively; Signal-to-noise ratio (SNR) of objectives after optimization are increased, it greatly enhances the robustness of the NCMB. The results demonstrate that the effectiveness of the methodology that based on Taguchi method and grey relational analysis.

**Keyword**：median barrier；multi-objective optimization; Taguchi；grey relational analysis；robustness

## 1 Introduction

Median barrier is a kind of energy-absorbing structure used to prevent cross-median crashes on two-lane or multi-lane highways. It works to contain and redirect the errant vehicle and relieve the occupant injury. Over the years, there were numerous hazardous events that the errant vehicles crossed median and collided with the oncoming vehicles in the opposite lane, which caused a lot of tragedies that killed many innocents and lost huge property. To mitigate the fatality and damage caused by the accidents, different types of barriers including concrete barrier, W-beam guardrail, and cable barrier have been developed. Among all the barrier types, concrete barrier behaves well in redirection, while it is in poor performance on reducing occupant injury. W-beam guardrail is the most widely used in highway, but it is especially expensive. Cable barrier is a road safety facility in rising use which is mild and less damaged. However, the deformation of this kind of barrier is very large, which may cause secondary crash events. The existing barriers installed on the U.S. highways are designed following the American Association of Transportation Officials (AASHTO) [1-3], and their performance is evaluated by the Manual for Assessing Safety Hardware (MASH) [4]. In China, the highway barriers are design referred to Design Specifications for Highway Safety Facilities (DSHSF, 2017) [5], and their performance is assessed by Standard for Safety Performance Evaluation of Highway Barriers (SSPEHB,2013) [6].

The researches on vehicle-barrier collision began in the 1920s. Full-scale crash test is the original method to investigate the barrier performance. Nowadays, it remains the most directive and effective means to study the crashworthiness of barrier. Bronstad and Burket [7] conducted six full-scale vehicle crash tests to evaluate the performance of timber weak-post system, it proved to be that the timber post satisfied the standard and the performance of timber posts was better than that of steel posts. In 1973, Hirsch et al [8] conducted full scale crash tests on three kinds of W-beam guardrail with different cross sections and lateral stiffness, and compared the performance of three kinds of guardrail. Ren and Vesenjak [9] tested the new road safety barrier according EU1317 standard by performing a full scale crash test. However, full-scale crash test is especially expensive and time-consuming for the test. Moreover, it's very difficult to do researches on multiple parameters. With the development of science and technology, computer simulation popularizes day by day. To reduce the cost of full-scale crash, nonlinear explicit finite element method is an advisable approach to evaluate the road safety barrier behavior. Since the 1970s, numerical simulation technology has been gradually applied to the numerical analysis in the field of automobile collision. In the 1980s, finite element method made great progress, and it rapidly became the main approach to simulate the vehicle-barrier collision process, which provided sufficient supplement and reference for the full-scale crash test.

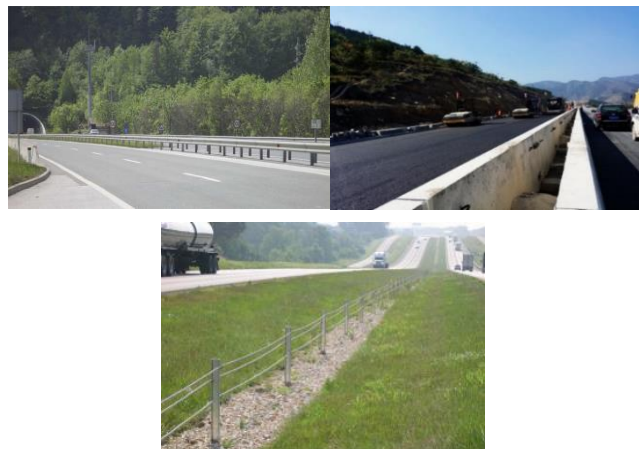

**Fig.1. Several common highway barriers:** （a）W-beam guardrail;（b）concrete barrier; (c) cable barrier

During the past decades, the finite element simulation commercial software such as LS-DYNA has been applied to design and evaluate the crashworthiness of the roadside hardware. Reid et al [10] used LS-DYNA to conducted a simulation study about guardrail, it's one of the earliest guardrail research associated LS-DYNA. After that, LS-DYNA is used to the roadside hardware crashworthiness investigation for its economy, efficiency and accuracy [11-17].

As to the design optimization problem involved in several and more parameters or multi-objective, finite element method is more prominent. Lei et al [18] used orthogonal experimental design method and finite element method to optimize the design parameters of discontinuous concrete guardrail. Hou et al [19] applied the RBF-MQ surrogate model to optimize the corrugated beam semi-rigid guardrails. Because of the complexity of the barrier system, there will be many abnormal phenomena such as tire snagging and vehicle rollover which does great damage to vehicles and occupants in the initial design of the structure. To avoid the tire snagging, Yin et al [20] used nonlinear finite element simulations combined with metamodeling-based design method to developed a $\eta$-shaped W-beam guardrail. Aimed at vehicle rollover happened on the bridge rail-to-guardrail transitions, Atahan and Cansiz [21] tackled the problem by increasing the W-beam height to 800mm. The traditional optimization method as above only takes the design variables into consideration, whereas it ignores the effect of the uncertainties on the performance of guardrail. In the design process, if the product is affected by the uncertainties, such as material properties fluctuates or test condition changes, the target response will be possible out of scope which can result in design failure. Robust design based on Taguchi method is proven available to lower the sensitivity of the product performance when influenced by the fluctuation of noise factor, and it is capable of improving the robustness [22-23].

It's easy to solve the problem of single objective optimization, but it doesn't work to tackle the multi-objective optimization problem. For this purpose, many researchers have conducted in-depth studies to deal with the multiple responses problem. Derringer and Suich [24] proposed a desirability function method, which converted the response value to the degree of desirability, and the total degree of multiple responses was realized by weighted geometric mean. He and Zhang [25] proposed an improved mahalanobis distance function method based on mahalanobis distance and quality lost function. However, in the field of multi-objective optimization, when a multi-objective optimization problem is transformed into a single-objective optimization problem, each objective function makes different contributions to the total objective function, so weight allocation is usually involved. What's more, the weight assignment problem is very difficult. Grey relational analysis provides a choice for dealing with complex multi - objective problems. Based on the grey correlation theory, the test data are treated with dimensionless corresponding quality characteristics, and then the grey relational coefficient between the two groups of data is calculated, and the relational degree is obtained through geometric average, so that the multi-objective problem is transformed into a single objective problem. Tosun [26] used grey relational analysis to determine the optimum parameters for multi-performance characteristics in drilling. Deepark and Davim [27] used grey relational method to reduce the surface taper and surface roughness, thus realized the multi-objective optimization of process paremeters in machining of graphite laced GFRP composite. Awale and Inamdar [28] optimized the turning parameters for lower cutting force, machining temperature, surface roughness and higher material removal rate simultaneously by using grey relational analysis. To the authors' best knowledge, many robustness designs only considered single objective, and most of the multi-objective optimization problems neglected the effect of uncertainties. Taguchi based on orthogonal is an efficient way to realize the robust optimization the discretization parameters of the structure, and grey relational analysis is application for the problem of multi-objective.

In this study, a new combined median barrier was developed. LS-DYNA based on the nonlinear finite element method was employed to evaluate the barrier performance. Considering the influence of uncertainties, combined with the Taguchi method and grey relational analysis, a robust multi-objective design model based on signal-to-noise ratio (SNR) and grey relational degree was established, and then the established model was used to realize the robust multi-objective optimization design of SBm level NCMB.

## 2 New combined median barrier modeling

### 2.1 Proposed of a new combined median barrier

A new combined median barrier (NCMB) was developed for reducing the severity of events that vehicles cross the median, which by enhancing the barrier foundation. It's composed of semi-rigid W-beam guardrail and rigid concrete structure. When vehicles impact the NCMB, the W-beam, crossbeam, and post of the semi-rigid guardrail as a resistance system to resist the vehicle collision. The concrete structure works not only as an enhanced barrier foundation, but also a structure for redirect the errant vehicle. In view of the drainage and the phenomenon of concrete expansion and contraction, the concrete structure is designed to discontinuous. The NCMB structure geometry model is shown in fig.2., where 1 denotes crossbeam, 2 denotes W-beam, 3 denotes concrete structure, 4 denotes post. The dimensions of the NCMB are show in Fig.3.

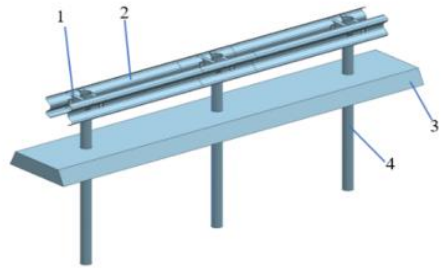

**Fig.2.** NCMB structure geometry model

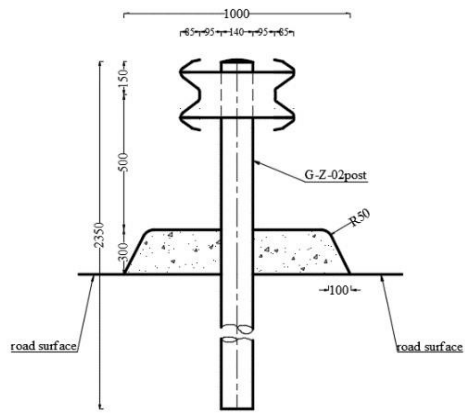

**Fig.3.** The dimensions of the NCMB

## 2.2 Finite element model of NCMB

Finite Element (FE) model of NCMB was developed in accordance with DSHSF, 2017. A detailed FE model of the NCMB is shown in Fig.4. In the FE model of barrier, commonly used elements include shell element, solid element and beam element. W-beam guardrail, crossbeam and post were simulated by shell elements, and the shell element used in this analysis is based on the Belytschko-Lin-Tsay shell formulation. Concrete structure was simulated by solid element, the beam element was used to simulate the reinforce.

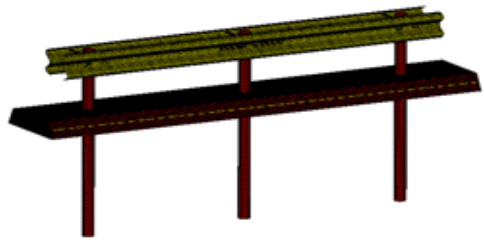

**Fig.4.** NCMB FE model

## 2.3 Finite element model of vehicles

In recent years, National Crash Analysis Center (NCAC) has developed many FE models for vehicles and road safety facilities. In this study, a FE model of Geo Metro Car developed and validated by NCAC was used to analysis. The complete vehicle model consists of 25037 elements and 28656 degrees of freedom (DOFs). The vehicle has a total mass of 900kg. The metal components and thin-walled structures of vehicles were simulated by shell elements; solid elements were applied for the modeling of power train, brake discs and callipers; beam elements were used to represent the rear suspension links [29]. According to the requirement of SSPEHB 2013, the mass of test vehicle is 1500kg, so the mass of car can reach the required quality by adding mass blocks. The Geo Metro Car FE model is shown in Fig.4.

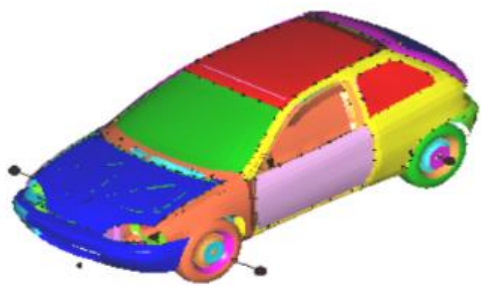

**Fig.4.** Car FE model

Medium bus FE model was developed at the Large-scale Structural Impact Laboratory of Changsha University of Science and Technology. It's reliability was validated by the collision between real vehicle and new cable barrier. The FE model has 207 parts, 53176 nodes and 53133 elements [30]. The bus FE model is shown in Fig.5, and the main technique parameters are shown in Tab.1.

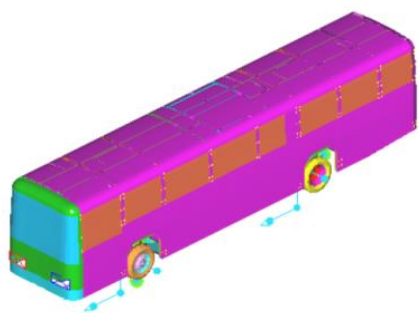

**Fig.5.** Medium bus FE model

**Tab.1.** Medium bus main technique parameters

| Parameters | Unit | Value |
|---|---|---|
| Total mass | kg | 10000 |
| Length×Width×Height | mm | 9020×2440×2880 |
| Wheelbase | mm | 1620/1820 |
| Centroid | mm | 1220 |

## 2.4 Coupling finite element model of vehicle and barrier

The coupling FE model of vehicle and barrier is shown in Figure 6. When the vehicle collides barrier in high-speed, it is usually accompanied by the occurrence of elastic-plastic material strain rate efficiency. This will have influence on the collision response. Plastic kinematic material can consider not only the strain rate, but also the failure of the materials. Therefore, plastic kinematic materials, namely MAT_PLASTIC_KINEMATIC constitutive model was used to simulate the post, W-beam, crossbeam and reinforce. The Cowper-Symonds model is capable of reacting to the strain rate effect, as described by Equation (1).

$$\sigma = \left[1 + \left(\frac{\dot{\varepsilon}}{C}\right)^{\frac{1}{p}}\right]\sigma_o \qquad (1)$$

Where $\sigma$ is the yield stress, $C$ and $P$ are the strain rate parameters for the Cowper-Symonds model, $\dot{\varepsilon}$ is the strain rate, $\sigma_0$ is the initial yield stress. For low carbon steel, the values of strain rate parameters are respectively C=40 and P=5[31].

Concrete structure was simulated by *MAT_JOHNSON_HOLMQUIST_ CONCRETE_TITLE model. Large strain, high strain rate and high pressure effect are considered in this model.

In the coupling FE model, the interaction between vehicle and barrier is complicated. The vehicle-barrier collision process is a complex contact process that involves geometry nonlinear, the effect of strain rate and size. In order to avoid the initial penetrations, there are 6 pairs of contact in the vehicle-barrier crash system simulation model, including the contact between vehicle and barrier, contact between vehicle and W-beam guardrail, contact between vehicle and concrete structure, contact for vehicle itself, contact for W-beam guardrail itself, contact for concrete structure itself. Tab.2. gives the contact definition of the crash system. SF represents the static friction coefficient, while DF represent the dynamic friction coefficient.

**Tab.2.** Vehicle-NCMB contact definition

| Contact elements | Contact order | SF | DF |
| --- | --- | --- | --- |
| Vehicle and W-beam guardrail | *CONTACT_AUTOMATIC_SURFACE_TO_SURFACE | 0.3 | 0.15 |
| Vehicle and concrete structure | *CONTACT_AUTOMATIC_SURFACE_TO_SURFACE | 0.3 | 0.15 |
| Vehicle itself | *CONTACT_AUTOMATIC_SINGLE_SURFACE | 0.3 | 0.15 |
| W-guardrail itself | *CONTACT_AUTOMATIC_SINGLE_SURFACE | 0.15 | 0.15 |
| Concrete itself | *CONTACT_AUTOMATIC_SINGLE_SURFACE | 0.15 | 0.15 |
| The whole barrier | *CONTACT_AUTOMATIC_SINGLE_SURFACE | 0.15 | 0.15 |

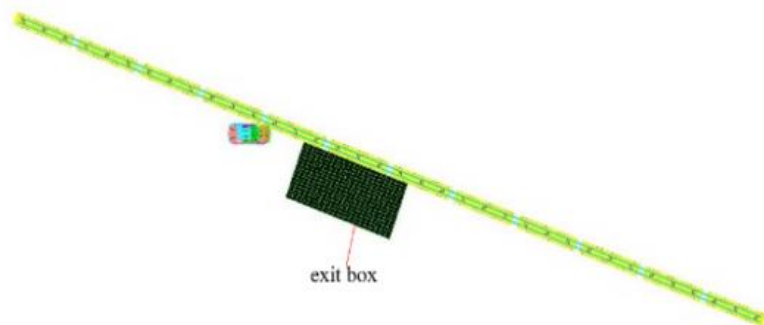

exit box

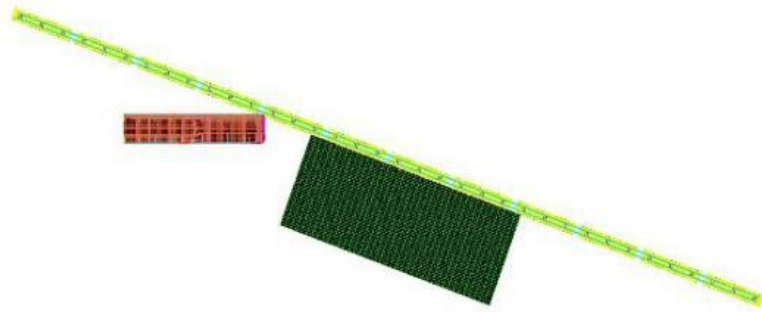

**Fig.6.** Vehicle-NCMB coupling FE model：(a) Car-NCMB FE model；(b) Medium bus-NCMB FE model

## 2.5 Impact conditions and Evaluation criteria
### 2.5.1 Impact conditions

FE simulations were used to investigate the safety performance of NCMB and realize the optimization of NCMB. In China, the design of highway barrier including roadside barrier and median barrier should meet the requirement specified by Standard for Safety Performance Evaluation of Highway Barriers(SSPEHB). According to the specification of SSPEHB, the FE crash simulations conditions as shown in Tab.3.

**Tab.3.** Test collision condition of the combined barrier

| Protection level | vehicle type | vehicle mass（t） | Impact velocity（km/h） | Impact angle（°） |
|---|---|---|---|---|
| SBm（240KJ） | Car | 1.5 | 100 | 20 |
| | Bus | 10 | 80 | 20 |

### 2.5.2 Evaluation criteria

The evaluation criteria for the vehicle responses after collision are also the SSPEHB. According to the SSPEHB, the safety performance of the NCMB is divided into qualitative evaluation criteria and quantitative evaluation criteria.

As for the qualitative evaluation criteria, it is mainly described by buffering performance and re-directive performance. Vehicles shall not roll over after collision, no turn around occurred, et al. Besides, the vehicle wheels after collision should remain in the exit box. The vehicle track after collision is shown in Figure 7.

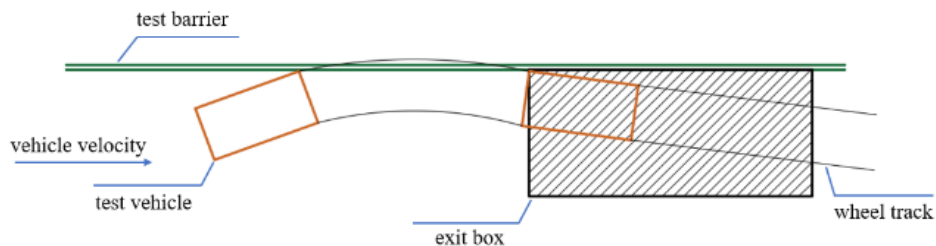

**Fig.7.** Track after collision

The quantitative evaluation mainly includes vehicle centroid acceleration(a) and maximum dynamic lateral defection of highway barriers(D). The A reflects the buffer performance of the barrier, it shall not be more than 20g, and the D is less than 1000mm.

# 3 Taguchi robust design method

## 3.1 Robust design theory

Product design process is susceptible to the influence of uncertainties, which can be generally classified into three categories: material properties uncertainties due to material heterogeneity, process differences and other factors; geometry size uncertainties caused by manufacturing error and assembly error; load case uncertainties caused by external environment. If there are slight changes in material properties and geometry size or load cases, the target response will fluctuate greatly, or even exceed the boundary constraints, which will have a destructive effect on product design [32].

Robust method is an efficient way that makes product or its performance insensitive to uncertainties. The uncertainties are called noise factors in the robust design, it's uncontrollable. In the product design stage, controllable factors and noise factors are considered. Controllable factors are also known as design variables, it can be controlled by human. The purpose of robust design is to find the optimal levels of controllable factors, which makes the product or production process have the strongest resistance to slight changes of noise factor. In the 1970s, Taguchi quality theory was proposed by Taguchi. Generally, the product quality in the robust design is described by mean quality loss function. It's main idea is to minimize the quality characteristics value variance from target value. The mean quality loss function is defined as Equation (2).

$$E\left[L(y)\right]=E\left[(y-T)^2\right] \qquad (2)$$

Where $y$ is the quality characteristics value, $L(y)$ is the quality loss function, $E[L(y)]$ is the mean quality loss function, T is the target value.

The above equation is transformed into equation (3) by adding $\bar{y}$ mean value, it minimizes not only the quality characteristics value variance from mean value, but also the mean value sensitivity from target value.

$$E\left[(y-T)^2\right] = E\left[(y-\bar{y})+(\bar{y}-T)\right]^2 \qquad (3)$$

Robust design schematic diagram is depicted as Fig.8.

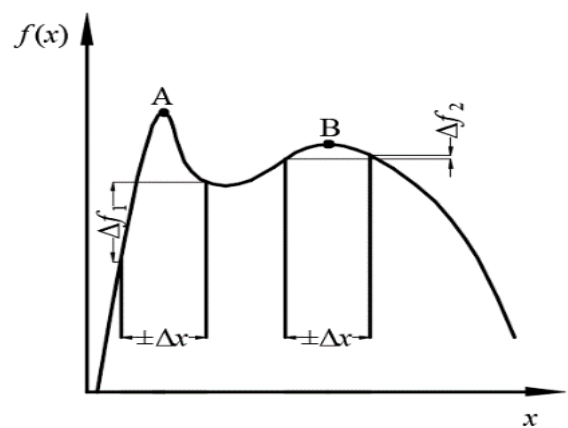

**Fig.8.** Robust design schematic diagram

## 3.2 Taguchi robust design method

Taguchi method includes three stages: system design, parameter design and tolerance design. Parameter design is the centre part of Taguchi method. It is an experimental method based on orthogonal design. Controllable factors and noise factors are arranged to orthogonal design respectively. Orthogonal array is composed of outside table and inside array. Noise factors are arranged to inside array, while controllable factors are arranged to outside array. Direct product array is generated by combining the outside array and inside array. In the cross point of controllable factor levels combination and noise factor levels combination will obtain an output response. In the end, the output response mean value is translated into signal-to-noise ratio (*SNR*). In general, the SNR refers to the ratio of effective information to noisy information in a signal system. In robust design, SNR is an index to measure the degree of robustness. The larger SNR value is, the more robust the product is. Conversely, the less robust it is. Taguchi robust parameter design array is shown in Table 4.

**Tab.4.** Parameter design array

| Table type | Inner table | | | | outside table | | | | | |
|---|---|---|---|---|---|---|---|---|---|---|
| Factor | Design variables | | | Sequence2 | | | | Noise factors | Response mean value $\overline{Y}$ | SNR |
| | | | | 1 | 2 | 3 | 4 | | | |
| | 1 | 2 | 3 | 1 | 2 | 2 | 1 | D | | |
| Sequence1 | | | | 1 | 2 | 1 | 2 | E | | |
| | A | B | C | 1 | 1 | 2 | 2 | F | | |
| 1 | 1 | 1 | 1 | $Y_{11}$ | $Y_{12}$ | $Y_{13}$ | $Y_{14}$ | | $\overline{Y}_1$ | $SNR_1$ |
| 2 | 1 | 2 | 3 | $Y_{21}$ | $Y_{22}$ | $Y_{23}$ | $Y_{24}$ | | $\overline{Y}_2$ | $SNR_2$ |
| 3 | 1 | 3 | 2 | $Y_{31}$ | $Y_{32}$ | $Y_{33}$ | $Y_{34}$ | | $\overline{Y}_3$ | $SNR_3$ |
| … | … | … | … | … | … | … | … | | … | … |

The quality characteristics of products are generally classified as Smaller the better, Larger the better, and Nominal is best. According to the characteristic difference of the quality characteristics, *SNR* is divided into these three different types. Assuming the output response is $y_1, y_2, \cdots, y_n$, T is the target value, three categories as given by following equations.

(1) Smaller the better:

$$SNR = -10\log\left(\frac{1}{n}\sum_{i=1}^{n} y_i^2\right) \quad (4)$$

(2) Larger the better:

$$SNR = -10\log\left(\frac{1}{n}\sum_{i=1}^{n}(\frac{1}{y_i})^2\right) \quad (5)$$

(3) Nominal is best:

$$SNR = -10\log\left(\frac{1}{n}\sum_{i=1}^{n}(y_i - T)^2\right) \quad (6)$$

# 4 Grey relational analysis

## 4.1 Grey relational theory

Grey relational analysis is an alternative method for multi-objective optimization. As to the multi-objective optimization problem, different target responses have different dimensions and orders of magnitude. Moreover, the objectives of multi-objective are opposite to each other, which makes the evaluation of multiple objectives unable to be uniform. It's also an efficient method for measuring the influence of the uncertain relations between each factor.

Because of the difference of range and unit of each objective, the original response data should be made being dimensionless.

When the objective is a characteristic of "Nominal is best", it can be normalized as Equation (7):

$$x_i^*(k) = \frac{\min_i\{E_{ik}, E_0\}}{\max_i\{E_{ik}, E_0\}} \tag{7}$$

If the objective is a characteristic of "Larger the better", it can be normalized as Equation (8):

$$x_i^*(k) = \frac{E_i(k) - \min_i E_i(k)}{\max_i E_i(k) - \min_i E_i(k)} \tag{8}$$

For the target value is a characteristic of "Smaller the better", the original response data should be normalized as following Equation (9):

$$x_i^*(k) = \frac{\max_i E_i(k) - E_i(k)}{\max_i E_i(k) - \min_i E_i(k)} \tag{9}$$

Where $i=1,\ldots, m$; $k=1,\ldots, n$. $m$ is the number of experiment, $n$ is the number of the objective. $E_i(k)$ is the original experiment response data, $\min_i E_i(k)$ is the minimal value of all the $i_{th}$ experiment response data, $\max_i E_i(k)$ is the maximum value of all the experiment response data. $E_0$ is the desired value. $x_i^*(k)$ is the sequence after being normalized.

After the original experiment response data is dimensionless, grey relational grade is calculated subsequently. Grey relational grade is a measurement of effects of each factor on the results. Grey relational grade is gained by applying grey relational theory as Equation (10).

$$\xi_i(k) = \frac{\min_i \min_k \Delta + \rho \min_i \min_k \Delta}{\Delta_0 + \rho \min_i \min_k \Delta} \tag{10}$$

Where $\xi_i(k)$ is the $\Delta$ is the deviation of the $i_{th}$ experiment response reference sequence and the $0_{th}$ comparative sequence, which can be expressed by equation: $\Delta = x_i^*(k) - x_0^*(k)$; $\rho$ is the distinguishing coefficient, it's defined in the range of 0-1.

$$\gamma_i = \frac{1}{n}\sum_1^n \xi_i \tag{11}$$

Where $\gamma_i$ exhibits the grey relational grade, $n$ exhibits the number of the objective.

## 4.2 Grey relational analysis based on Signal-to-Noise ratio

Combining the *SNR* and grey relational theory, a new multi-objective robust design model is building. Firstly, the original response of quality characteristic is transformed into *SNR*. Apparently, *SNR* is a characteristic of larger the better. So the *SNR* will be normalized further by the grey relational analysis which is a characteristic of larger the better. The new multi-objective robust design model as following

equations:

Firstly, the initial dimensionless of *SNR* is shown as Equation (12).

$$\eta_i(k) = \frac{SNR_i(k) - \min_i SNR_i(k)}{\max_i SNR_i(k) - \min_i SNR_i(k)} \quad (12)$$

Where $i=1,…, m$; $k=1,…, n$, $m$ is the number of test, $n$ is the number of objective, $SNR_i$ is the test data, $\min_i SNR_i$ is the minimum value of $SNR_i$, $\max_i SNR_i$ is the maximum of $SNR_i$, $SNR_0$ is the target value。$\eta_i(k)$ is the level after dimensionless。

Secondary, calculation of grey relational coefficient is shown as Equation (13).

$$\xi^{(\eta)}_i(k) = \frac{\min_i \min_k \Delta_\eta + \rho \max_i \max_k \Delta_\eta}{\Delta_\eta + \rho \max_i \max_k \Delta_\eta} \quad (13)$$

Where $\xi^{(\eta)}_i(k)$ is the grey relational coefficient; $\Delta_\eta$ is the *SNR* difference of referential sequence and comparative sequence,: $\Delta_\eta = \eta_i(k) - \eta_0(k)$; $\rho$ is the identification coefficient limited in the range of 0~1。

Lastly, grey relational grade is calculated by Equation (14).

$$\gamma^{(\eta)}_i = \frac{1}{n} \sum_1^n \xi^{(\eta)}_i \quad (14)$$

Where $\gamma_i$ is the grey relational grade, $n$ is the target number.

# 5 Results and discussion

## 5.1 Definition of optimization objectives

The robust optimization objectives are the vehicle centroid acceleration(a) and maximum dynamic lateral defection of highway barriers(D). The lower the value of the A is, the less damage to the occupant. It is of characteristic of "Smaller the better". At the same time, the smaller the D value is, the better the barrier performance is. The smaller D value can prevent the deformable components intruding into the opposite lane. It is of characteristic of "Larger the better".

## 5.2 Definition of design variables and noise factors
### 5.2.1 Definition of design variables

Post thickness, cross beam thickness, corrugated beam thickness and post spacing were chosen to be the design variables preliminarily. We performed a univariate sensitivity analysis for these four variables. The factors selected and levels are shown in the Table 5.

Table 5 Factors and levels of univariate sensitivity analysis                                    mm

| Levels | Factors | | | |
| --- | --- | --- | --- | --- |
| | Post thickness A | Cross beam thickness B | Corrugated beam thickness C | Post spacing D |
| 1 | 3.5 | 3.5 | 3 | 1500 |
| 2 | 4.5 | 4.5 | 4 | 2000 |
| 3 | 5.5 | 5.5 | 5 | 2500 |

Changing the post thickness and the other factors remained the same, post thickness was performed a univariate sensitivity analysis. Results are shown in the Figure 9.

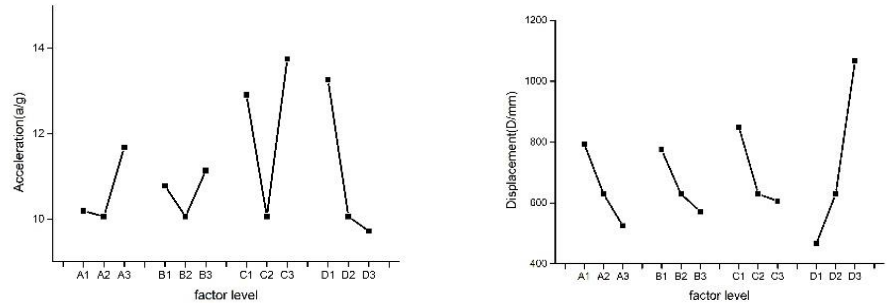

**Fig.9.** Univariate sensitivity analysis：(a) univariate sensitivity analysis of acceleration；(b) univariate sensitivity analysis of deflection

According to the univariate sensitivity analysis above, Post thickness, corrugated beam thickness and post spacing were determined to be the design variables of robust optimization design.

**5.2.2 Definition of noise factors**

In the process of vehicle impacts the NCMB, the velocity and mass of the vehicle is unable to control, so the velocity ($v$) and mass ($m$) were as the noise factors of the robust design. The material characteristics such as yield stress are affected by the equipment accuracy and the material defects, and it's uncontrollable, so the yield stress was as another noise factor. The noise factors and levels are shown in the Table 6.

**Tab.6.** Noise factors in different levels

| Levels | Factors | | |
|---|---|---|---|
|  | Veloctity($v$/m/s) | Mass($m$/kg) | Yield Stress($\sigma$/MPa) |
| 1 | 80 | 10000 | 235 |
| 2 | 84 | 10300 | 240 |

## 5.3 Analysis of the ignal-to-noise ratios

### 5.3.1 Parameter design based on Taguchi method

Orthogonal experimental design was applied to parameter design based on Taguchi method. 2 orthogonal tables were required for the parameter design, one called inner table, it was used to arrange the design variables. While the other is outside table, it is used for noise factors. According to the levels of design variables and noise factors, $L_9(3^3)$ orthogonal table was selected for design variables, and $L_4(2^3)$ orthogonal table was selected for noise factors, and a total of $9\times4=36$ experiments were carried out.

**Tab.7.** Experiment arrangement for Parameter design based on Taguchi method

| Table type | Inner table | | | outside table | | | | | | |
|---|---|---|---|---|---|---|---|---|---|---|
| Factor | Design variables | | | Sequence2 | | | | Noise factors | Response mean value $\overline{Y}$ | SNR |
|  |  |  |  | 1 | 2 | 3 | 4 |  |  |  |
| Sequence1 | 1 | 2 | 3 | 1 | 2 | 2 | 1 | D |  |  |
|  |  |  |  | 1 | 2 | 1 | 2 | E |  |  |

| | A | B | C | 1 | 1 | 2 | 2 | F | | |
|---|---|---|---|---|---|---|---|---|---|---|
| 1 | 1 | 1 | 1 | $Y_{11}$ | $Y_{12}$ | $Y_{13}$ | $Y_{14}$ | | $\overline{Y}_1$ | $SNR_1$ |
| 2 | 1 | 2 | 3 | $Y_{21}$ | $Y_{22}$ | $Y_{23}$ | $Y_{24}$ | | $\overline{Y}_2$ | $SNR_2$ |
| 3 | 1 | 3 | 2 | $Y_{31}$ | $Y_{32}$ | $Y_{33}$ | $Y_{34}$ | | $\overline{Y}_3$ | $SNR_3$ |
| 4 | 2 | 1 | 2 | $Y_{41}$ | $Y_{42}$ | $Y_{43}$ | $Y_{44}$ | | $\overline{Y}_4$ | $SNR_4$ |
| 5 | 2 | 2 | 1 | $Y_{51}$ | $Y_{52}$ | $Y_{53}$ | $Y_{54}$ | | $\overline{Y}_5$ | $SNR_5$ |
| 6 | 2 | 3 | 3 | $Y_{61}$ | $Y_{62}$ | $Y_{63}$ | $Y_{64}$ | | $\overline{Y}_6$ | $SNR_6$ |
| 7 | 3 | 1 | 3 | $Y_{71}$ | $Y_{72}$ | $Y_{73}$ | $Y_{74}$ | | $\overline{Y}_7$ | $SNR_7$ |
| 8 | 3 | 2 | 2 | $Y_{81}$ | $Y_{82}$ | $Y_{83}$ | $Y_{84}$ | | $\overline{Y}_8$ | $SNR_8$ |
| 9 | 3 | 3 | 1 | $Y_{91}$ | $Y_{92}$ | $Y_{93}$ | $Y_{94}$ | | $\overline{Y}_9$ | $SNR_9$ |

According to the experiment design, the response of acceleration is shown in the Table 8. As the acceleration is of characteristics of "smaller the better", the corresponding *SNR* was calculated by Equation (4).

**Tab.8.** Response of acceleration

| Sequence | Acceleration(a/g) | | | | | |
|---|---|---|---|---|---|---|
| | $a_1$ | $a_2$ | $a_3$ | $a_4$ | $a$ | SNR |
| 1 | 16.13 | 17.93 | 16.09 | 19.58 | 17.43 | -24.86 |
| 2 | 11.75 | 7.99 | 11.36 | 8.39 | 11.61 | -20.01 |
| 3 | 11.92 | 14.29 | 12.18 | 15.33 | 13.43 | -22.61 |
| 4 | 14.35 | 9.43 | 7.48 | 4.30 | 8.89 | -19.65 |
| 5 | 9.03 | 11.51 | 10.23 | 9.27 | 10.01 | -20.05 |
| 6 | 9.55 | 11.79 | 10.75 | 8.76 | 10.21 | -20.24 |
| 7 | 10.57 | 6.84 | 7.99 | 9.49 | 8.72 | -18.92 |
| 8 | 10.88 | 8.17 | 10.97 | 11.20 | 10.31 | -20.32 |
| 9 | 5.93 | 4.64 | 4.85 | 4.41 | 4.96 | -13.96 |

The response of maximum dynamic lateral defection of highway barriers is shown in Table 9. The "smaller the better" Equation (4) was also chosen for the calculation of maximum dynamic lateral deflection.

**Tab.9.** Response of maximum dynamic lateral defection

| Sequence | Deflection(D/mm) | | | | | |
|---|---|---|---|---|---|---|
| | $D_1$ | $D_2$ | $D_3$ | $D_4$ | $D$ | S/N |
| 1 | 806 | 892 | 838 | 826 | 840.5 | -58.50 |
| 2 | 719 | 945 | 1071 | 794 | 882.3 | -59.01 |
| 3 | 662 | 964 | 702 | 760 | 772.0 | -57.85 |
| 4 | 1017 | 829 | 743 | 696 | 821.3 | -58.39 |
| 5 | 478 | 734 | 623 | 512 | 586.8 | -55.49 |
| 6 | 569 | 727 | 1009 | 669 | 743.5 | -57.63 |
| 7 | 731 | 728 | 704 | 649 | 703.0 | -56.95 |

| 8 | 545 | 658 | 649 | 345 | 549.3 | -55.02 |
| 9 | 356 | 373 | 364 | 353 | 361.5 | -51.16 |

## 5.4 Grey relational analysis

Grey relational analysis is applied for analysis of ambiguity of a system and deficiency of information. The *SNR* was taken as "larger is better ", so the obtained *SNR* were further normalized with Equation (12). The results are shown in Table 10.

**Tab.10.** Normalized *SNR* of the series

| Serial No. | Acceleration | Deflection |
| --- | --- | --- |
| 1 | 0 | 0.076 |
| 2 | 0.445 | 0 |
| 3 | 0.206 | 0.157 |
| 4 | 0.478 | 0.089 |
| 5 | 0.441 | 0.455 |
| 6 | 0.424 | 0.185 |
| 7 | 0.545 | 0.301 |
| 8 | 0.417 | 0.514 |
| 9 | 1 | 1 |

In order to measure the importance of each parameter on the experimental results, the grey relational coefficients (GRC) of each parameter should be obtained by Equation (13). To further analysis, grey relational grade (GRD) were computed by Equation (14). Grey relational coefficients and grey relational grade values of the series is shown in Table 11.

**Tab.11.** Grey relational coefficients and grey relational grade values of the series

| Serial No. | GRC | | GRD | Order |
| --- | --- | --- | --- | --- |
| | Acceleration | Deflection | | |
| 1 | 0.333 | 0.351 | 0.3420 | 9 |
| 2 | 0.474 | 0.333 | 0.4035 | 7 |
| 3 | 0.386 | 0.372 | 0.3790 | 8 |
| 4 | 0.489 | 0.354 | 0.4215 | 6 |
| 5 | 0.472 | 0.478 | 0.4750 | 3 |
| 6 | 0.456 | 0.380 | 0.4225 | 5 |
| 7 | 0.524 | 0.417 | 0.4705 | 4 |
| 8 | 0.462 | 0.507 | 0.4845 | 2 |
| 9 | 1 | 1 | 1 | 1 |

## 5.5 Analysis of optimized results

### 5.5.1 Analysis of range

The range reflects the maximum dispersion of a set of data. The greater the range, the greater the influence of this factor on the result. It can be calculated by Equation (15).

$$R = \max\{\overline{K_i}\} - \min\{\overline{K_i}\} \qquad (15)$$

Where $K_i$ is the mean deviation in level $i$.

The result of range analysis is shown in Tab 12. It can be seen that the significance of each factor by order on the NCMB performance are: A>B>C. And the optimal combination is: $A_3B_3C_1$.

**Tab.12.** Results of analysis of range

| Factors | Levels | | | Range |
|---|---|---|---|---|
| | 1 | 2 | 3 | |
| A | 0.3478 | 0.4397 | 0.6517 | 0.2769 |
| B | 0.4113 | 0.4543 | 0.6005 | 0.1892 |
| C | 0.6057 | 0.4338 | 0.4322 | 0.1735 |

**5.5.2 Analysis of variance**

ANOVA has been used for obtaining the contribution ratio of each parameter on the target function. Moreover, another aim of using ANOVA method is to verify the importance level of the parameters obtained from optimization. It can be performed by Equation (16)-Equation (18). The result of ANOVA is shown in Table 13.

$$\mu = \frac{1}{mn}\sum_{i=1}^{m}\sum_{j=1}^{n}x_{ij} \qquad (16)$$

$$S_T = \sum_{i=1}^{m}\sum_{j=1}^{n}(x_{ij} - \mu)^2 \qquad (17)$$

$$F = \frac{S_1/f_1}{S_2/f_2} \qquad (18)$$

**Table 13** Results of ANOVA

| Factors | Deviation sum of squares | df | Mean square | F | Significance |
|---|---|---|---|---|---|
| A | 0.078 | 2 | 0.039 | 0.839 | 0.544 |
| B | 0.080 | 2 | 0.040 | 0.854 | 0.539 |
| C | 0.061 | 2 | 0.030 | 0.653 | 0.605 |

From the results of ANOVA, we know that the contribution ratio of design variables on the GRG by order: C>A>B. Thus, the final design scheme can be determined as: $A_3B_3C_1$.

**5.5.3 Analysis of the robust design results**

In order to verify the accuracy of multi-objective robust design results, a comparison analysis will be performed between the optimized variables and the original variables. The velocity of vehicle, mass of vehicle, yield stress of NCMB components is taken as the design variables, and taking 2 levels for each variable, a L4(3$^2$) orthogonal experiment design is performed for the NCMB before and after optimization.

**Tab.14.** Comparison analysis of before and after optimization

| Serials. No | Factors Levels | | | Acceleration | | Deflection | |
|---|---|---|---|---|---|---|---|
| | A | B | C | Before | After | Before | After |
| 1 | 1 | 1 | 1 | 11.79 | 5.93 | 742 | 356 |
| 2 | 2 | 1 | 2 | 12.30 | 4.85 | 692 | 364 |
| 3 | 1 | 2 | 2 | 10.58 | 4.41 | 741 | 353 |
| 4 | 1 | 2 | 1 | 9.79 | 4.64 | 664 | 373 |
| Mean | | | | 11.115 | 4.96 | 709.75 | 361.5 |
| Standard deviation | | | | 0.999 | 0.526 | 33.259 | 7.762 |
| *SNR* | | | | -20.92 | -13.91 | -57.02 | -51.16 |

The results show that the mean of acceleration and deflection after optimization is 4.96g and 361.5mm respectively, which were reduced by 55.4% and 49.1%; the standard deviation of acceleration and deflection after optimization is 0.526 and 7.762 respectively, which were reduced by 47.3% and 76.7% respectively; the *SNR* of acceleration and deflection after optimization is -13.91 and -51.16, which were increased by 33.5% and 10.3%.

As shown in Fig.10., the vehicle after collision drives out follow the exit box, the redirection of the vehicle is successful. It meets the requirement of SSPEHB.

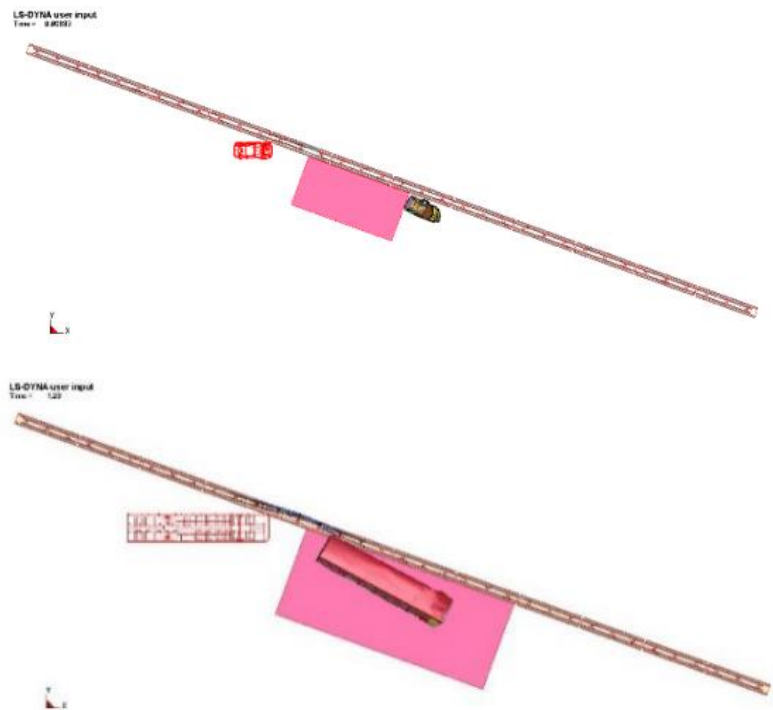

**Fig.10.** Redirection of the vehicle:(a) redirection of car; (b) redirection of medium bus

# 6 Conclusion

The multi-objective optimization of new combined median barrier was achieved by means of Taguchi and grey relational analysis. The conclusions can be drawn based on the experiment and response analysis:

(1) The post thickness, corrugated beam thickness and post spacing have great effects on the crashworthiness of NCMB.

(2) Generally speaking, the greater the stiffness of the barrier, the worse buffering to the occupant. However, in this study, with the increase of the thickness of the post and corrugated beam, the acceleration does not fall but increase. It indicates that the stiffness is not always positively correlated with occupant damage. Moreover, the deflection of the barrier is not always negatively correlated with the acceleration.

(3) The results after robust design based on Taguchi method and grey relational analysis, the acceleration and deflection are greatly reduced, the crashworthiness of the barrier is improved; the mean value and standard deviation value of the response are decreased, and the $SNR$ is increased. It indicates that the robustness of NCMB structure has been improved.

## Acknowledgments

This work is supported by National Key Research and Development Program of the 13th Five-Year Plan of China(No.2016YFB1200602-33) and the National Natural Science Foundation of China(No.52072048).